\documentstyle[12pt]{article}

\newenvironment{references}[1]{\begin{footnotesize}%
\begin{thebibliography}{#1}}
{\end{thebibliography}\end{footnotesize}
}
\newtheorem{theorem}{Theorem}[section]

\def\z{{\bf z}}
\def\rho{{\varrho}}

\def\mm{{\cal {\bf M}}}

\title{ Note on universal algorithms for learning theory}

\author{
 Karol  Dziedziul,\\
Barbara Wolnik,\\
}

\begin{document}

\maketitle


\abstract{ We propose the general way of study the universal estimator for the
regression problem in learning theory considered in
\cite{BCDDT1} and
\cite{BCDDT2}.
This new approch allows us, for example, to improve the results from
\cite{BCDDT1}.}

{\bf Keywords:} {nonparametric regression, learning theory.}

Mathematics Subject Classification:  68T05, 41A36, 41A45, 62G05.


\section{Introduction}

 I am recalling this paper since it is somehow important. We wrote this  paper
  to improve one of the results from \cite{BCDDT1}. It was written in \cite{BCDDT1} that
" It has been communicated to us
by Lucien Birg\'e that one can derive from one of his forthcoming papers (Birg\'e, 2004) [published in (2006)] that for any class $\Theta$ satisfying (5), (namely $\Theta\subset L^2(X,\rho_X)$ with a condition on the entropy
number which is similar to the assumption (10)), there is an estimator $f_z$ satisfying
\[
E \|f_\rho-f_{\z}\|^2=O \left({1\over m}\right)^{2s\over 1+2s}.
\]
whenever $f_\rho \in \Theta$". Our paper shows how to construct an estimator  with straightforward reasoning and I can not find the mentioned  proof so I recall the proof of the result.

\noindent   S. Cucker and S. Smale in their paper \cite{CS}
determined the scope of  the learning theory. We would like to
present a general approach which corresponds to the papers \cite{BCDDT1} and
\cite{BCDDT2}. The problem is the following. Let $X=[0,1]^d$ and
$Y=[-A,A]$. On a product space $Z=X\times Y$ it is  unknown
probability Borel measure $\rho$. We shall assume that the
marginal probability measure $\rho_X(S)=\rho(S\times Y)$ on $X$ is
a Borel measure. We have
\[
d\rho(x,y)=d\rho(y|x) d\rho_X(x).
\]

We are given the data $\z \subset Z$ of $m$ independent random
observation $z_j=(x_j,y_j)$, $j=1,2,\ldots, m$ identically
distributed according to $\rho$. We are interested in estimating
{\bf the regression function}
\[
f_\rho(x):=\int_Y yd\rho(y|x)
\]
in $L^2(X,\rho_X)$ norm which will be denoted by $\|\cdot\|$.

\bigskip

\noindent
To do it let $\mm =\{ M_v\}_{v\in T}$ denote any family of measurable
functions on $X$ such that for all $v\in T$
\begin{equation}\label{01}
0\leq M_v(x)\leq 1,\;\;\; x\in X
\end{equation}
and
\begin{equation}\label{suma}
\sum_{v\in T} M_v(x)=1,\;\;\; x\in X.
\end{equation}
One of examples of this kind of family $\mm$ is the family $\{\chi_I\}_{I\in
T}$, where $\chi_I$ denotes the indicator function of $I$ and $\{ I:\;
I\in T\}$ is any partition of the set $X$ (in
\cite{BCDDT1} the sets $I$ are dyadic cubes). Another example we get if we
considere a triangulation of $X$ with the vertices
$\{ v\}_{v\in T}$. To define piecewise linear and continuous function
corresponding to every vertex $v\in T$ it is sufficient to define
such function on vertices. We define basis
\[
M_v(w)=\cases{1& for vertices $w= v$\cr 0 & for $w\neq v$.}
\]
It is not hard to check that family $\{ M_v\}_{v\in T}$ satisfies (\ref{01})
and (\ref{suma}).

\bigskip

 Now for a given family $\mm$ we define the operator
\[
Q_{\mm}f(x)=\sum_{v\in T}c_v(f) M_v(x),
\]
where
\[
c_v(f)={\alpha_v(f)\over \rho_v}, \quad \alpha_v(f)=\int_X f M_v d\rho_X,
\quad \rho_v=\int_X  M_v d\rho_X
\]
and the estimator
\[
f_{\z}(x)=\sum_{v\in T}c_v(\z) M_v(x),
\]
where
\[
 c_v(\z)={\alpha_v(\z)\over
\rho_v(\z)},
\]
\[
 \alpha_v(\z)={1\over
m}\sum_{j=1}^m y_j M_v(x_j), \quad \rho_v(\z)={1\over
m}\sum_{j=1}^m M_v(x_j).
\]
If $\rho_v=0$ then we define $c_v=0$ and if
$\rho_v(\z)=0$ then we put $c_v(\z)=0.$  Note also that
$E\alpha_v(\z)=\alpha_v$ (here and subsequently
$\alpha_v:=\alpha_v(f_{\rho})$, $c_v:=c_v(f_{\rho})$) and
$E\rho_v(\z)=\rho_v$. Moreover
\[
Var(yM_v(x))\leq \int_Z y^2 M^2_v(x)d\rho(x,y) \leq A^2 \int_X
M_v^2(x) d\rho_X(x),
\]
hence
\begin{equation}\label{1}
Var(yM_v(x)) \leq A^2 \int_X M_v(x) d\rho_X(x)=A^2 \rho_v,
\end{equation}
\begin{equation}\label{2}
Var(M_v(x))\leq E(M_v(x))^2 \leq  E(M_v(x))=\rho_v.
\end{equation}

\noindent
Therefore by Bernstein's inequality, see for instance 
\cite{BCDDT1}  we have for any $\epsilon >0$
\begin{equation}\label{Ba}
Prob\left\{ \left|{\alpha_v}-{\alpha_v(\z)}\right|
\geq \epsilon \right\}\leq 2e^{-\frac{3m\epsilon^2}{6A^2\rho_v+4A\epsilon}},
\end{equation}
\begin{equation}\label{Br}
Prob\left\{ \left|{\rho_v}-{\rho_v(\z)}\right|\geq {\epsilon} \right\}\leq
2e^{-\frac{3m\epsilon^2}{6\rho_v+2\epsilon}}.
\end{equation}

\bigskip

The main result of this paper is
\begin{theorem}
For any family $\mm$
\[
E \|Q_{\mm}f_\rho-f_{\z}\|^2=O({N\over m}),
\]
 where $N=|T|$.
\end{theorem}

\noindent
The new idea of the proof presented below allows us to improve the result from
\cite{BCDDT1} (in Corollary 2.2 \cite{BCDDT1}
the above expectation is bounded by $O(\frac{N}{m}\cdot \log N)$).

\noindent
{\sc Proof}. By (\ref{01}), (\ref{suma}) and the convexity of
the square functions we have
\[
E\|Q_{\mm} f_\rho-f_{\z}\|^2 \leq \int_X \sum_{v\in T}E|c_v-c_v(\z)|^2
M_v(x) d\rho_X(x)
\]
\[
= \sum_{v\in T}E|c_v-c_v(\z)|^2 \rho_v.
\]
Note that  if $\rho_v=0$ then $E\rho_v(\z)=0$ hence $\rho_v(\z)=0$
$\rho^m$ a.e. Consequently
\[
E\|Q_{\mm} f_\rho-f_{\z}\|^2  \leq \sum_{v\in T,
\rho_v>0}E|c_v-c_v(\z)|^2 \rho_v.
\]
Let us fix $v$ such that $\rho_v>0$. We can write
\[
E|c_v-c_v(\z)|^2=\int_{\rho_v(\z)>0}|c_v-c_v(\z)|^2
+\int_{\rho_v(\z)=0}|c_v|^2.
\]
Note that if $\rho_v(\z)=0$ $\rho^m$ a.e. then for all $j$
$M_v(x_j)=0$, hence $\alpha_v(\z)=0$ $\rho^m$ a.e. Thus
\[
E|c_v-c_v(\z)|^2=\int_{\rho_v(\z)>0}|c_v-c_v(\z)|^2
+\int_{\rho_v(\z)=0}|{\alpha_v-\alpha_v(\z)\over \rho_v}|^2.
\]
For $b\neq 0$ and $t\neq 0$ we use the simple inequality
\begin{equation}\label{n}
\left|\frac{a}{b}-\frac{s}{t}\right|\leq\frac{1}{|b|}|a-s|+\frac{|s|}{|bt|}|t-b|
\end{equation}
to get
\begin{equation}\label{33}
\left|{a\over b}-{s\over t} \right|^2\leq 2 {|a-s|^2\over
b^2}+2{1\over b^2}{s^2\over t^2}|t-b|^2,
\end{equation}
which gives in particular that
\[
\left| \frac{a_v}{\rho_v}-\frac{a_v(\z)}{\rho_v(\z)}\right|^2\leq
2\frac{|a_v-a_v(\z)|^2}{\rho_v^2}+2\left(\frac{a_v(\z)}{\rho_v(\z)}\right)^2
\frac{|\rho_v-\rho_v(\z)|^2}{\rho_v^2}.
\]
For $\rho_v(\z)>0$ we have
\[
{\alpha_v(\z)^2\over \rho_v(\z)^2}\leq A^2.
\]
thus
\[
E|c_v-c_v(\z)|^2\leq
{3\over m\rho_v^2} Var
(yM_v(x))+
{2A^2\over m\rho_v^2}Var(M_v(x))).
\]
Consequently
\[
E\|Q_T f_\rho-f_{\z}\|^2 \leq C \sum_{v\in T} {1\over m\rho_v^2}
(Var (yM_v(x))+Var(M_v(x)))\rho_v.
\]
By (\ref{1}) and (\ref{2}) we get
\[
E\|Q_T f_\rho-f_{\z}\|^2 \leq O( \sum_{v\in T} {1\over
m})=O({N\over m})
\]
and this finishes the proof.

\vskip 1cm

Let us note that if we take $N=m^{{1\over 1+2s}}$ for fixed
$s>0$ then
\begin{equation}\label{3}
E \|Q_{\mm}f_\rho-f_{\z}\|^2=O \left({1\over m}\right)^{2s\over
1+2s}.
\end{equation}
To unify approach to linear and nonlinear approach in estimation
let us introduce the sets ${\cal A}^s$ similar to definition given
in \cite{BCDDT1}. We have that $f\in {\cal A}^s$, $0<s$ (in fact
it makes sense to consider $0<s\leq 2$) if $f\in L^2(\rho_X)$ and
there is $C$ such that for all $N$ there is a family $\mm=\{ M_v\}_{v\in T}$
with properties (\ref{01}) and (\ref{suma}) such that $N=|T|$ and
\begin{equation}\label{4}
\|f-Q_{\mm}f\|\leq C N^{-s}.
\end{equation}

\noindent
By Theorem 1.2, (\ref{3}) and
(\ref{4}) and since
\[
E\|f_\rho-f_{\z}\|^2\leq
2E\|f_\rho-Q_{\mm}f_{\rho}\|^2+2E\|Q_{\mm}f_{\rho}-f_{\z}\|^2
\]
we get the optimal rate of estimation (see \cite{DKPT}). This
approach  improves the rate of estimation in  (\cite{BCDDT1}).

\begin{theorem}
Let $f_\rho\in {\cal A}^s$ and let $\mm$ be the family from
definition of space ${\cal A}^s$ such that $N=|T|=[m^{{1\over
1+2s}}]$. Then
\[
E \|f_\rho-f_{\z}\|^2=O \left({1\over m}\right)^{2s\over 1+2s}.
\]
\end{theorem}

Finally, we will show the general version of the Theorem 2.1 in \cite{BCDDT1}.
Our proof is very analogous but partially simplified, so we present
it for the sake of completeness.
We improve the constant in estimation.

\begin{theorem}
For any family $\mm$ and any $\eta>0$
\begin{equation}\label{o1}
Prob\{ \|Q_{\mm}f_\rho-f_{\z}\|>\eta\}\leq 4N e^{-cm\eta^2\over N},
\end{equation}
where $N:=|T|$ and $c$ depends only on $A$.

\end{theorem}

{\sc Proof}. By the convexity of the square function we have that
\begin{equation}\label{s}
\|Q_{\mm} f_\rho-f_{\z}\|^2 \leq \int_X \sum_{v\in T}|c_v-c_v(\z)|^2
M_v(x) d\rho_X(x)
= \sum_{v\in T}|c_v-c_v(\z)|^2 \rho_v.
\end{equation}

\noindent
This gives
\[
Prob\{\|Q_{\mm} f_\rho-f_{\z}\| >\eta\} \leq
Prob\{ \sum_{v\in T}|c_v-c_v(\z)|^2 \rho_v>\eta^2\}
\]
\[
\leq\sum_{v\in T}Prob\{|c_v-c_v(\z)|>\frac{\eta}{\sqrt{N\rho_v}}\}.
\]

\noindent
Let us note that
\[
Prob\{|c_v-c_v(\z)|>\frac{\eta}{\sqrt{N\rho_v}}\}=0
\]
provided $\rho_v\leq\frac{\eta^2}{4A^2N}$. To see this it is enough to
transform this assumption to the form $\frac{\eta}{\sqrt{N\rho_v}}\geq 2A$ and
recall that $|c_v|$ and $|c_v(\z)|$ are less than $A$.

\noindent
Therefore we can write
\[
Prob\{\|Q_{\mm} f_\rho-f_{\z}\| >\eta\} \leq
\sum_{v:\;\;\rho_v>\frac{\eta^2}{4A^2N}}Prob\{|c_v-c_v(\z)|>\frac{\eta}{\sqrt{N\rho_v}}\}.
\]

\noindent
To estimate the last sum let us note that if
\[
|\alpha_v(\z)-\alpha_v|\leq \frac{\rho_v\eta}{4\sqrt{N\rho_v}}
\]
and
\[
|\rho_v(\z)-\rho_v|\leq \frac{\rho_v\eta}{4A\sqrt{N\rho_v}}
\]
then ( we know that $\rho_v>\frac{\eta^2}{4A^2N}$)
\[
|\rho_v(\z)-\rho_v|\leq \frac{\rho_v\eta}{4A\sqrt{N
\frac{\eta^2}{4A^2N}
}}=\frac{1}{2}\rho_v
\]
(this gives in particular that $|\rho_v(\z)|\geq \frac{1}{2}\rho_v$) and using
(\ref{n})
we get
\[
|c_v(\z)-c_v|=\left|\frac{\alpha_v(\z)}{\rho_v(\z)}-\frac{\alpha_v}{\rho_v}\right|
\]
\[
\leq
\frac{1}{|\rho_v(\z)|}|\alpha_v(\z)-\alpha_v|+\frac{|\alpha_v|}{|\rho_v(\z)|\rho_v}|\rho_v(\z)-\rho_v|
\]
\[
\leq \frac{1}{\frac{1}{2}\rho_v}\cdot
\frac{\rho_v\eta}{4\sqrt{N\rho_v}}+
\frac{A}{\frac{1}{2}\rho_v}\cdot
\frac{\rho_v\eta}{4A\sqrt{N\rho_v}}=
\frac{\eta}{\sqrt{N\rho_v}}.
\]

\noindent
Therefore
\[
Prob\left\{|c_v-c_v(\z)|>\frac{\eta}{\sqrt{N\rho_v}}\right\}
\]
\[
\leq
Prob\left\{|\alpha_v(\z)-\alpha_v|> \frac{\rho_v\eta}{4\sqrt{N\rho_v}}\right\}+
Prob\left\{|\rho_v(\z)-\rho_v|> \frac{\rho_v\eta}{4A\sqrt{N\rho_v}}\right\}.
\]
If we first use (\ref{Ba}), (\ref{Br}) and then the fact that $\frac{\eta}{\sqrt{N\rho_v}}\leq
2A$ we get finally
\[
Prob\{\|Q_{\mm} f_\rho-f_{\z}\| >\eta\} \leq
\]
\[
\leq\sum_{v:\;\;\rho_v>\frac{\eta^2}{4A^2N}}
\left(
2e^{-\frac{3m\eta^2}{16N(6A^2+A\frac{\eta}{\sqrt{N\rho_v}})}}
+2e^{-\frac{3m\eta^2}{16A^2N(6+\frac{1}{2A}\cdot\frac{\eta}{\sqrt{N\rho_v}})}}
\right)
\]
\[
\leq
\sum_{v:\;\;\rho_v>\frac{\eta^2}{4A^2N}}
2\left(
e^{-\frac{3}{128}\cdot\frac{m\eta^2}{NA^2}}
+e^{-\frac{3}{112}\cdot\frac{m\eta^2}{NA^2}}
\right)
\leq 4Ne^{-\frac{3}{128A^2}\cdot\frac{m\eta^2}{N}}.
\]
which complete the proof of (\ref{o1}) with $c=\frac{3}{128A^2}$.

Karol  Dziedziul,\\
karol.dziedziul@pg.edu.pl\\

Barbara Wolnik,\\
barbara.wolnik@mat.ug.edu.pl

\end{document}